\documentclass{article}

\usepackage{hyperref}
\usepackage{cleveref}
\hypersetup{hypertexnames = false, bookmarksdepth = 2, bookmarksopen = true, colorlinks, linkcolor = black, citecolor = black, urlcolor = black, pdfstartview={XYZ null null 1}}

\usepackage{amsfonts}
\usepackage{amsthm}
\usepackage{amssymb}
\usepackage{enumitem}
\usepackage{mathtools}
\usepackage{thmtools}
\usepackage{tikz}
\usepackage{tikz-cd}
\usepackage[textsize=small]{todonotes}
\usepackage{xparse}
\usepackage{xspace}

\def\Ascr{{\mathcal A}}

\def\Fscr{{\mathcal F}}

\def\Oscr{{\mathcal O}}

\def\Rscr{{\mathcal R}}
\def\Sscr{{\mathcal S}}
\def\Tscr{{\mathcal T}}

\theoremstyle{definition}
\declaretheorem[name=Theorem]{theorem}

\newtheorem{definition}[theorem]{Definition}
\newtheorem{example}[theorem]{Example}

\newtheorem{lemma}[theorem]{Lemma}

\newtheorem{proposition}[theorem]{Proposition}
\newtheorem{remark}[theorem]{Remark}

\mathchardef\mhyphen="2D
\newcommand\dash{\nobreakdash-\hspace{0pt}}

\newcommand\bounded{\ensuremath{\mathrm{b}}}
\newcommand\braid{\ensuremath{\mathrm{B}}}
\newcommand\Sigmabraid{\ensuremath{\Sigma\mathrm{B}}}
\newcommand\LLL{\ensuremath{\mathbf{L}}}
\newcommand\RRR{\ensuremath{\mathbf{R}}}

\DeclareMathOperator\Aut{Aut}

\DeclareMathOperator\Bl{Bl}
\DeclareMathOperator\BS{BS}
\DeclareMathOperator\codim{codim}
\DeclareMathOperator\coh{coh}

\DeclareMathOperator\derived{\mathbf{D}}
\DeclareMathOperator\Ext{Ext}
\DeclareMathOperator\GKdim{GK\,dim}
\DeclareMathOperator\gldim{gl\,dim}
\DeclareMathOperator\Gr{Gr}
\DeclareMathOperator\gr{gr}
\DeclareMathOperator\Hom{Hom}
\DeclareMathOperator\id{id}
\DeclareMathOperator{\Ker}{ker}
\DeclareMathOperator\Kzero{K_0}
\DeclareMathOperator\leftmutate{L}
\DeclareMathOperator\Mat{Mat}

\DeclareMathOperator\Pic{Pic}
\DeclareMathOperator\Proj{Proj}
\DeclareMathOperator\Qcoh{Qcoh}
\DeclareMathOperator\QGr{QGr}
\DeclareMathOperator\qgr{qgr}
\DeclareMathOperator\relSpec{\underline{Spec}}

\DeclareMathOperator\Sym{Sym}
\DeclareMathOperator\rk{rk}

\DeclareMathOperator\ZZ{Z}

\title{Construction of noncommutative surfaces with exceptional collections of length 4}
\author{Pieter Belmans \and Dennis Presotto}

\begin{document}

\maketitle

\begin{abstract}
  Recently de Thanhoffer de V\"olcsey and Van den Bergh classified the Euler forms on a free abelian group of rank~4 having the properties of the Euler form of a smooth projective surface.
  There are two types of solutions: one corresponding to~$\mathbb{P}^1\times\mathbb{P}^1$ (and noncommutative quadrics), and an infinite family indexed by the natural numbers. For~$m=0,1$ there are commutative and noncommutative surfaces having this Euler form, whilst for~$m\geq 2$ there are no commutative surfaces. In this paper we construct sheaves of maximal orders on surfaces having these Euler forms, giving a geometric construction for their numerical blowups.

\end{abstract}

\tableofcontents


\section{Introduction}
In a recent paper de Thanhoffer de V\"olcsey and Van den Bergh provide a \emph{numerical} classification of possibly noncommutative surfaces with an exceptional sequence of length~4~\cite{1607.04246v1}. Their classification describes the possible bilinear forms on a free abelian group of rank~4 mimicking the properties of the numerical Grothendieck group and Euler form on a smooth projective surface.

Those properties are described as follows: for a finitely generated free abelian group~$\Lambda$, a nondegenerate bilinear form~$\langle-,-\rangle\colon\Lambda\times\Lambda\to\mathbb{Z}$ and an automorphism~$s\in\Aut(\Lambda)$ we will ask that
\begin{description}[labelwidth=\widthof{\bfseries Serre automorphism\quad},align=right]
  \item[Serre automorphism] $\langle x,s(y)\rangle=\langle y,x \rangle$ for~$x,y \in \Lambda$;
  \item[unipotency] $s-\id_\Lambda$ is nilpotent;
  \item[rank] $\rk(s-\id_\Lambda)=2$;
\end{description}
By the nondegeneracy we know that if we choose a basis for~$\Lambda$, and express the bilinear form (resp.~the Serre automorphism) as the Cartan or Gram matrix~$M$ (resp.~the Coxeter matrix~$C$), we have the relation
\begin{equation}
  C=-M^{-1}M^{\mathrm{t}},
\end{equation}
so it suffices to specify the bilinear form.

The reason for considering these properties is explained in~\cite{1607.04246v1}: it can be shown that the action of the Serre functor on the Grothendieck group for all smooth projective surfaces has the extra property that~$s-\id_\Lambda$ has rank precisely~2, whilst the unipotency of the Serre functor holds in complete generality~\cite[lemma~3.1]{MR1230966}. Independently, Kuznetsov developed a similar notion in \cite{MR3691718}.

Before giving the classification, recall that mutation and shifting of exceptional collections gives an action of the signed braid group~$\Sigmabraid_n$, and we will only be interested in the classification (of bilinear forms) up to this action.

For rank~3 the analogous problem is classical and is described by the Markov equation. In that case the only solution is given by~$\mathbb{P}^2$, and its noncommutative analogues. The structure of the numerical Grothendieck group in this case can be read off from \cref{example:commutative-beilinson,example:noncommutative-beilinson}.

For rank~4 there are more solutions, which are described by de Thanhoffer de V\"olcsey and Van den Bergh in their main theorem~\cite[theorem~A]{1607.04246v1}. They show that, up to the action of the signed braid group, the matrices
\begin{equation}
  \tag{A}
  \label{equation:quadric-matrix}
  \begin{pmatrix}
    1 & 2 & 2 & 4 \\
    0 & 1 & 0 & 2 \\
    0 & 0 & 1 & 2 \\
    0 & 0 & 0 & 1
  \end{pmatrix}
\end{equation}
and
\begin{equation}
  \tag{$\mathrm{B}_m$}
  \label{equation:family-matrix}
  \begin{pmatrix}
    1 & m & 2m & m \\
    0 & 1 & 3 & 3 \\
    0 & 0 & 1 & 3 \\
    0 & 0 & 0 & 1
  \end{pmatrix}
\end{equation}
for~$m\in\mathbb{N}$, describe all possible bilinear forms~$\langle-,-\rangle$ on~$\mathbb{Z}^4$ satisfying the properties.

The case~\eqref{equation:quadric-matrix} corresponds to the quadric surface~$\mathbb{P}^1\times\mathbb{P}^1$ (and its noncommutative analogues~\cite{MR2836401}). Using \cref{proposition:family-P2-matrix} the cases \eqref{equation:family-matrix} with~$m=0$ (resp.~$1$) correspond to the disjoint union of~$\mathbb{P}^2$ (and its noncommutative analogues~\cite{MR1086882,MR1230966}) with a point, resp.~the blowup of~$\mathbb{P}^2$ in a point (and its noncommutative analogues~\cite{MR1846352}).

For~\eqref{equation:family-matrix} with~$m=2$ de Thanhoffer de V\"olcsey--Presotto constructed families of noncommutative~$\mathbb{P}^1$\dash bundles on~$\mathbb{P}^1$ of rank~$(1,4)$, and showed that these give the correct Euler form~\cite{1503.03992v4}.

In this paper we give a streamlined construction of families of noncommutative surfaces with Grothendieck group~$\mathbb{Z}^4$ and Euler form~\eqref{equation:family-matrix} for all~$m\geq 2$ using completely different methods. This is achieved by constructing a sheaf of maximal orders on~$\Bl_x\mathbb{P}^2$, as done in \cref{section:construction}. The main technique here is a noncommutative generalisation of Orlov's blowup formula, as given in \cref{theorem:orlov-for-orders}, and which is probably of independent interest.

The result of the constructions in this paper can be summarised as follows.
\begin{theorem}
  \label{theorem:main-theorem-summary}
  For every~$m\geq 2$ there exist maximal orders on~$\Bl_x\mathbb{P}^2\cong\mathbb{F}_1$ whose Euler form is of type \eqref{equation:family-matrix}.
\end{theorem}
In particular, we provide an actual geometric construction for the numerical blowups of \cite{1607.04246v1}. There are three degrees of freedom in our construction, which is the expected number of degrees of freedom, but a complete classification of noncommutative surfaces of rank~4 is out of reach.

In \cref{section:order-properties} we give some properties of the orders we have constructed in the context of the minimal model program for orders. Especially for~$m=2$ they turn out to be interesting, as we get interesting new examples of so called \emph{half ruled orders}. These are also the only ones which are del Pezzo, which is somewhat unexpected as the analogous construction for blowing up a point on the ramification locus always gives a del Pezzo order.

For the case of~$m=2$ there are now two constructions of an abelian category with the prescribed properties: we have that it is possible to view the abelian category both as a blowup, and a~$\mathbb{P}^1$\dash bundle, just like it is possible for~$m=1$ to view~$\Bl_x\mathbb{P}^2$ as the Hirzebruch surface~$\mathbb{F}_1\coloneqq\mathbb{P}(\mathcal{O}_{\mathbb{P}^1}\oplus\mathcal{O}_{\mathbb{P}^1}(-1))$ and vice versa. In a future work~\cite{belmans-presotto-vandenbergh-comparison} we give a comparison of the two constructions.

\paragraph{Acknowledgements}
The authors would like to thank Daniel Chan and Colin Ingalls for interesting discussions regarding the properties of our orders from the point of view of the minimal model program. The authors would like to thank Alexander Kuznetsov and Colin Ingalls for interesting comments on a draft version of this paper.

Both authors were supported by an FWO PhD fellowship.

{\ }

{\small This is the accepted version of the following article: Construction of non‐commutative surfaces with exceptional collections of length 4, which has been published in final form at Journal of London Mathematical Society.}

\section{Preliminaries}
\label{section:preliminaries}

\subsection{Artin--Schelter regular algebras}
Artin and Schelter introduced in~\cite{MR917738} a class of graded algebras that were to serve as the noncommutative analogues of the polynomial ring~$k[x_1, \ldots ,x_d]$. They are defined as follows.
\begin{definition}
  Let~$A$ be a connected graded~$k$\dash algebra. Then we say~$A$ is \emph{Artin--Schelter regular} of dimension~$d$ if
  \begin{enumerate}
    \item $\gldim A = d$;
    \item $\GKdim A = d$;
    \item $A$ is Gorenstein (with respect to the integer~$d$), i.e.~there exists~$l\in\mathbb{Z}$ such that
      \begin{equation}
        \Ext_{\Gr A}^i(k_A,A)\cong
        \begin{cases}
          {}_Ak(l) & i=d, \\
          0 & i\neq d.
        \end{cases}
      \end{equation}
  \end{enumerate}
\end{definition}

Of particular importance in the study of noncommutative surfaces is the case~$d=3$, and we will restrict ourselves to this case. Moreover we only consider 3-dimensional AS-regular algebras which are generated in degree~1. For these algebras, it turns out there are two possible Hilbert series~\cite[theorem~1.5(i)]{MR917738}, one of which is precisely that of~$k[x,y,z]$. These algebras are referred to as quadratic AS-regular algebras and the associated abelian category~$\qgr A$ is called a \emph{noncommutative plane}. The other class of algebras is related to~$\mathbb{P}^1\times\mathbb{P}^1$, and is of no role here.

Quadratic AS-regular algebras were classified in terms of triples of geometric data~\cite[definition~4.5]{MR1086882}.
\begin{definition}
  An \emph{elliptic triple} is a triple~$(C,\sigma,\mathcal{L})$ where
  \begin{enumerate}
    \item $C$ is a divisor of degree~3\ in~$\mathbb{P}^2$;
    \item $\sigma\in\Aut(C)$;
    \item $\mathcal{L}$ is a very ample line bundle of degree~3 on~$C$.
  \end{enumerate}
  We say that it is \emph{regular} if moreover
  \begin{equation}
    \mathcal{L}\otimes\left( \sigma^*\circ\sigma^*(\mathcal{L}) \right)\cong\sigma^*(\mathcal{L})\otimes\sigma^*(\mathcal{L}).
  \end{equation}
\end{definition}

The classification of noncommutative planes using regular triples is originally due to Artin--Tate--Van den Bergh~\cite[\S4]{MR1086882}, and later done using different techniques by Bondal--Polishchuk~\cite{MR1230966}.

\begin{example}
  \label{example:sklyanin}
  The generic case in the classification is given by the \emph{Sklyanin algebra}. The geometric data in this case is an elliptic curve, and the automorphism is given by a translation.

  It is well known that a Sklyanin algebra can be written as the quotient of~$k\langle x,y,z\rangle$ by the ideal generated by
  \begin{equation}
    \label{equation:sklyanin-equations}
    \left\{
      \begin{aligned}
        axy + byx + cz^2 &= 0 \\
        ayz + bzy + cx^2 &= 0 \\
        azx + bxz + cy^2 &= 0
      \end{aligned}
    \right.
  \end{equation}
  where~$[a:b:c] \in \mathbb{P}^2 \setminus S$ and~$S$ is an explicitly known finite set of~12~points.

\end{example}

We will be particularly interested in the case where the Sklyanin algebra~$A$ is finite over its center. This property is visible in the geometric data~\cite[theorem~7.1]{MR1128218}.
\begin{theorem}[Artin--Tate--Van den Bergh]
  \label{theorem:finite-over-center}
  The algebra~$A$ is finite over its center if and only if the automorphism in the associated elliptic triple is of finite order.
\end{theorem}

\begin{remark}
  \label{remark:finite-order-automorphisms}
  By inspecting the automorphism groups of cubic curves it can be seen that there are only~4~point schemes (out of~9~point schemes appearing in \cite[table 1]{MR1230966}) that are allowed in order for the automorphism to be of finite order. These are the elliptic curves, the nodal cubic, the triangle of lines and a conic and line in general position. For an alternative approach using the classification of maximal orders, see \cref{remark:artin-mumford-sequence}.
\end{remark}

We will use the following properties of the algebras which are finite over their center:
\begin{enumerate}
  \item we can obtain the category~$\qgr A$ as the category~$\coh\mathcal{S}$, where~$\mathcal{S}$ is a sheaf of maximal orders on~$\mathbb{P}^2$, as explained in \cref{lemma:equivalence-qgr-coh};
  \item there exist so called fat points in~$\qgr A$, as explained in \cref{subsection:fat-point-modules}.
\end{enumerate}
These will allow us to construct a sheaf of maximal orders on~$\mathbb{F}_1\cong\Bl_x\mathbb{P}^2$ with the desired properties in \cref{section:construction}.

\subsection{Semiorthogonal decompositions}
\label{subsection:semiorthogonal-decompositions}
We will also use the notion of semiorthogonal decomposition and exceptional sequences. We will denote~$\mathcal{T}$ a~$k$\dash linear triangulated category, and all constructions are~$k$\dash linear. For more details one is referred to~\cite{MR3545926}.

\begin{definition}
  A \emph{semiorthogonal decomposition} of~$\mathcal{T}$ is a sequence~$(\mathcal{S}_1,\dotsc,\mathcal{S}_n)$ of full triangulated subcategories of~$\mathcal{T}$, such that there exists a filtration
  \begin{equation}
    0=\mathcal{T}_0\subseteq\mathcal{T}_1\subseteq\dotsc\subseteq\mathcal{T}_n=\mathcal{T}
  \end{equation}
  where~$\mathcal{T}_i\subseteq\mathcal{T}$ is a (left) admissible subcategory such that~$\mathcal{S}_i\cong\mathcal{T}_i/\mathcal{T}_{i-1}$. We will denote this by
  \begin{equation}
    \label{equation:sod}
    \mathcal{T}=\left\langle \mathcal{S}_1,\dotsc,\mathcal{S}_n \right\rangle.
  \end{equation}
\end{definition}

A special case of a semiorthogonal decomposition is provided by a full exceptional sequence. In this case we have that the~$\mathcal{S}_i$'s are equivalent to the derived category of the base field~$k$.

\begin{definition}
  We say that~$E\in\mathcal{T}$ is an \emph{exceptional object} if
  \begin{equation}
    \Hom_{\mathcal{T}}(E,E[m])\cong
    \begin{cases}
      k & m=0 \\
      0 & m\neq 0
    \end{cases}.
  \end{equation}
  A sequence~$(E_1,\dotsc,E_n)$ of objects~$E_i\in\mathcal{T}$ is an \emph{exceptional sequence} if each~$E_i$ is exceptional, and~$\Hom_{\mathcal{T}}(E_i,E_j[m])=0$ for all~$i>j$ and~$m$. It is said to be \emph{strong} if moreover~$\Hom_{\mathcal{T}}(E_i,E_j[m])=0$ for all~$m\neq 0$. It is said to be \emph{full} if it generates the category~$\mathcal{T}$, or equivalently
  \begin{equation}
    \mathcal{T}=\left\langle \langle E_1\rangle,\dotsc,\langle E_n\rangle \right\rangle.
  \end{equation}
  is a semiorthogonal decomposition.
\end{definition}

The main property of semiorthogonal decompositions that we will use in this paper is that they are sent to direct sums by so called additive invariants. In particular for the Grothendieck group we get in the situation of \eqref{equation:sod} that
\begin{equation}
  \Kzero(\mathcal{T})\cong\bigoplus_{i=1}^n\Kzero(\mathcal{S}_i).
\end{equation}
In particular, if~$\mathcal{T}$ has a full exceptional collection of length~$n$, then
\begin{equation}
  \Kzero(\mathcal{T})\cong\mathbb{Z}^{\oplus n}.
\end{equation}

The main example of a full and strong exceptional collection, and also the example that motivated the construction in \cref{section:construction} is Beilinson's collection on~$\mathbb{P}^2$~\cite{MR509388}.
\begin{example}
  \label{example:commutative-beilinson}
  The derived category of~$\mathbb{P}^2$ has a full and strong exceptional collection
  \begin{equation}
    \derived^\bounded(\mathbb{P}^2)=
    \left\langle
      \mathcal{O}_{\mathbb{P}^2},
      \mathcal{O}_{\mathbb{P}^2}(1),
      \mathcal{O}_{\mathbb{P}^2}(2)
    \right\rangle
  \end{equation}
  whose quiver is
  \begin{equation}
    \label{eq:beilinson-quiver}
    \begin{tikzpicture}[scale = 2, baseline = (current bounding box.center)]
      \draw[circle, minimum size = 7pt, inner sep = 0pt]
        (0,0) node (v0) {}
        (1,0) node (v1) {}
        (2,0) node (v2) {};

      \draw (v0) circle (1pt)  
            (v1) circle (1pt)  
            (v2) circle (1pt); 

      \draw[->] (v0) edge [bend left = 45]  node[fill = white] {$x_0$} (v1)
                (v0) edge                   node[fill = white] {$y_0$} (v1)
                (v0) edge [bend left = -45] node [fill = white] {$z_0$} (v1)
                (v1) edge [bend left = 45]  node [fill = white] {$x_1$} (v2)
                (v1) edge                   node [fill = white] {$y_1$} (v2)
                (v1) edge [bend left = -45] node [fill = white] {$z_1$} (v2);
    \end{tikzpicture},
  \end{equation}
  and the relations are
  \begin{equation}
    \label{eq:relations-P2}
    \left\{
      \begin{aligned}
        x_0y_1&=y_0x_1 \\
        x_0z_1&=z_0x_1 \\
        y_0z_1&=z_0y_1.
      \end{aligned}
    \right.
  \end{equation}
  In particular, this means that~$\Kzero(\derived^\bounded(\mathbb{P}^2))\cong\mathbb{Z}^{\oplus3}$, and we can read off that the Gram (or Cartan) matrix is
  \begin{equation}
    \label{equation:cartan-P2}
    M=
    \begin{pmatrix}
      1 & 3 & 6 \\
      0 & 1 & 3 \\
      0 & 0 & 1
    \end{pmatrix},
  \end{equation}
  whilst the Coxeter matrix is
  \begin{equation}
    \label{equation:coxeter-P2}
    C=
    \begin{pmatrix}
      -10 & -6 & -3 \\
      15 & 8 & 3 \\
      -6 & -3 & -1
    \end{pmatrix}.
  \end{equation}
\end{example}

This example can be generalised to noncommutative~$\mathbb{P}^2$'s in the following way.
\begin{example}
  \label{example:noncommutative-beilinson}
  The derived category of~$\qgr A$, where~$A$ is a quadratic~3\dash dimensional Artin--Schelter regular algebra, has a well-known full and strong exceptional collection mimicking that of Beilinson for~$\mathbb{P}^2$, given by
  \begin{equation}
    \derived^\bounded(\qgr A)=
    \left\langle
      \pi A,
      \pi A(1),
      \pi A(2)
    \right\rangle
  \end{equation}
  where~$\pi$ denotes the quotient functor~$\pi\colon\gr A\to\qgr A$, and~$A(i)$ denotes the grading shift of~$A$. More details can be found in~\cite[theorem~7.1]{1411.7770v1}.

  The quiver has the same shape as \eqref{eq:beilinson-quiver}, and the relations can be read off from the presentation of~$A$ as a quotient of~$k\langle x,y,z\rangle$ by~3~quadratic relations as in~\eqref{eq:relations-P2}, see also \cref{lemma:structure-DbS}.

  For instance in the case of the Sklyanin algebra of \cref{example:sklyanin} they are
  \begin{equation}
    \label{equation:sklyanin-exceptional-collection-equations}
    \left\{
      \begin{aligned}
        ax_0y_1 + by_0x_1 + cz_0z_1 &= 0 \\
        ay_0z_1 + bz_0y_1 + cx_0x_1 &= 0 \\
        az_0x_1 + bx_0z_1 + cy_0y_1 &= 0
      \end{aligned}
    \right.
  \end{equation}
  The Cartan and Coxeter matrices~$A$ and~$C$ describing the Euler form and the Serre functor only depend on the structure of the quiver with relations, and this stays the same, so we obtain the matrices from~\eqref{equation:cartan-P2} and~\eqref{equation:coxeter-P2}. In \cite{1607.04246v1} it is explained how up to the signed braid group action introduced in \cref{subsection:mutation} this is the only solution of rank~3.
\end{example}

\subsection{Mutation}
\label{subsection:mutation}
We quickly recall the theory of mutations of exceptional sequences.


\begin{definition}
  \label{definition:mutation}
  Let~$\Tscr$ be an~$\Ext$-finite triangulated category and let~$(E,F)$ be an exceptional pair of objects in~$\Tscr$. We define the \emph{left mutation}~$\leftmutate_E F$ as the cone of the morphism (since the pair is exceptional, the cone is unique up to \emph{unique} isomorphism).
  \begin{equation}
    \label{equation:mutation}
    \Hom_\Tscr(E,F)\otimes E \rightarrow F \rightarrow \leftmutate_E F
  \end{equation}
  If~$\mathbb{E}\coloneqq(E_1, \ldots E_n)$ is an exceptional collection in~$\Tscr$ we define the \emph{mutation at~$i$} to be the exceptional collection~$(E_1, \ldots, \leftmutate_{E_i}E_{i+1}, E_i, \ldots, E_n)$.

  These mutations can be interpreted as an action of the braid group on~$n$~strings, denoted~$\braid_n$, on the set of all exceptional collections. To see this, let~$\sigma_1, \ldots, \sigma_{n-1}$ be the standard generators for~$\braid_n$, then~$\sigma_i$ acts on an exceptional collection via mutation at~$i$, i.e.
  \begin{equation}
    \sigma_i(\mathbb{E}):= (E_1, \ldots,\leftmutate_{E_i}E_{i+1}, E_i, \ldots, E_n).
  \end{equation}
\end{definition}


\begin{remark}
  By a celebrated theorem by Kuleshov and Orlov~\cite{MR1286839} we know that for a del Pezzo surface~$X$ the braid group~$\braid_m$ (where~$m = \rk\Kzero(X)$) acts transitively on the set of exceptional collections in~$\derived^\bounded(X)$.
\end{remark}

Inspired by~\cite{1607.04246v1} we also consider the action of the \emph{signed} braid group, which also takes shifting into account.

\begin{definition}
  The \emph{signed braid group}~$\Sigmabraid_n$ is the semidirect product~$\braid_n\rtimes(\mathbb{Z}/2\mathbb{Z})^n$, where~$(\mathbb{Z}/2\mathbb{Z})^n$ acts on~$\braid_n$ by considering the quotient~$\braid_n\twoheadrightarrow\Sym_n$. As such, the signed braid group has~$2n-1$ generators:
  \begin{itemize}
    \item $n-1$ generators~$\sigma_1,\ldots,\sigma_{n-1}$, as for the braid group~$\braid_n$,
    \item $n$ generators~$\epsilon_1,\ldots,\epsilon_n$, as for~$(\mathbb{Z}/2\mathbb{Z})^n$.
  \end{itemize}
  These generators satisfy the relations
  \begin{equation}
    \left\{
      \begin{aligned}
        \sigma_i\sigma_j
        &=\sigma_j\sigma_i
        & |i-j|\geq 2 \\
        \sigma_i\sigma_{i+1}\sigma_i
        &=\sigma_{i+1}\sigma_i\sigma_{i+1}
        & i=1,\dotsc,n-2 \\
        \epsilon_i^2
        &=1 \\
        \epsilon_i\epsilon_j
        &=\epsilon_j\epsilon_i \\
        \epsilon_i\sigma_i\epsilon_{i+1}
        &=\sigma_i
        & i=1,\dotsc,n-1.
      \end{aligned}
    \right.
  \end{equation}
\end{definition}

In~\cite[\S4]{1607.04246v1} the rules for computing the action of~$\Sigmabraid_n$ on a bilinear form~$\langle-,-\rangle$ are given. These rules naturally generalize the induced action of~$\braid_n$ on the Euler form on the Grothendieck group~$\Kzero(\Tscr)$. We will construct an abelian category of ``geometric origin'' which has a Grothendieck group that is in the same orbit as the matrix~\eqref{equation:family-matrix}. To do so we will use the representative of the equivalence class found in the next proposition.

\begin{proposition}
  \label{proposition:family-P2-matrix}
  The matrices \eqref{equation:family-matrix} are mutation equivalent to the matrices
  \begin{equation}
    \tag{$\mathrm{B}'_m$}
    \label{equation:family-P2-matrix}
    \begin{pmatrix}
      1 & 3 & 6 & m \\
      0 & 1 & 3 & m \\
      0 & 0 & 1 & m \\
      0 & 0 & 0 & 1
    \end{pmatrix}
  \end{equation}
  for~$n\in\mathbb{N}$.

  \begin{proof}
    We have that~$\epsilon_1\epsilon_3\sigma_3\sigma_1\sigma_2\sigma_3$ sends~\eqref{equation:family-P2-matrix} to~\eqref{equation:family-matrix}. The mutations~$\sigma_1\sigma_2\sigma_3$ provide a shift in the helix, whilst the mutation~$\sigma_3$ at that point corresponds to the mutation that sends~$(\mathcal{O}_{\mathbb{P}^2},\mathcal{O}_{\mathbb{P}^2}(1),\mathcal{O}_{\mathbb{P}^2}(2))$ to~$(\mathcal{O}_{\mathbb{P}^2},\mathrm{T}_{\mathbb{P}^2}(-1),\mathcal{O}_{\mathbb{P}^2}(1))$.

    The intermediate steps are
    \begin{equation}
      \begin{aligned}
        \sigma_3\eqref{equation:family-P2-matrix}
        &=
        \begin{pmatrix}
          1 & 3 & -5m & 6 \\
          0 & 1 & -2m & 3 \\
          0 & 0 & 1 & -m \\
          0 & 0 & 0 & 1
        \end{pmatrix}
        \\
        \sigma_2\sigma_3\eqref{equation:family-P2-matrix}
        &=
        \begin{pmatrix}
          1 & m & 3 & 6 \\
          0 & 1 & 2m & 5m \\
          0 & 0 & 1 & 3 \\
          0 & 0 & 0 & 1
        \end{pmatrix}
        \\
        \sigma_1\sigma_2\sigma_3\eqref{equation:family-P2-matrix}
        &=
        \begin{pmatrix}
          1 & -m & -m & -m \\
          0 & 1 & 3 & 6 \\
          0 & 0 & 1 & 3 \\
          0 & 0 & 0 & 1
        \end{pmatrix}
        \\
        \sigma_3\sigma_1\sigma_2\sigma_3\eqref{equation:family-P2-matrix}
        &=
        \begin{pmatrix}
          1 & -m & -2m & -m \\
          0 & 1 & -3 & 3 \\
          0 & 0 & 1 & -3 \\
          0 & 0 & 0 & 1
        \end{pmatrix}
        \\
        \epsilon_3\sigma_3\sigma_1\sigma_2\sigma_3\eqref{equation:family-P2-matrix}
        &=
        \begin{pmatrix}
          1 & -m & -2m & -m \\
          0 & 1 & 3 & 3 \\
          0 & 0 & 1 & 3 \\
          0 & 0 & 0 & 1
        \end{pmatrix}
      \end{aligned}
    \end{equation}
  \end{proof}
\end{proposition}

\subsection{Fat point modules}
\label{subsection:fat-point-modules}
According to~\cite[\S7]{MR1273836} we take:
\begin{definition}
  \label{definition:fat-point-module}
  A \emph{fat point module} for~$A$ is a graded module~$F$ satisfying the following properties:
  \begin{enumerate}
    \item\label{enumerate:fat-point-module-1} $F$ is generated by~$F_0$
    \item\label{enumerate:fat-point-module-2} The Hilbert function~$\dim_k F_n$ is a constant~$\geq 2$, which is called the \emph{multiplicity}.
    \item\label{enumerate:fat-point-module-3} $F$ has no nonzero finite-dimensional submodules.
    \item\label{enumerate:fat-point-module-4} $\pi F \in \qgr A$ is simple
  \end{enumerate}
  A \emph{fat point} is the isomorphism class of a fat point module in~$\qgr A$.
\end{definition}

The following result tells us that fat point modules for quadratic Artin--Schelter regular algebras which are finite over their center behave particularly well.

By the classification of those algebras which are finite over their center we know that such an algebra is described by an elliptic triple~$(E,\sigma,\mathcal{L})$ where~$\sigma$ is an automorphism of finite order. We will denote
\begin{equation}
  s\coloneqq\min\{k\mid\sigma^{k,*}(\mathcal{L})\cong\mathcal{L}\}.
\end{equation}
It is this integer, and not the order of~$\sigma$ (which will be denoted~$n$) that is the important invariant of the triple~$(E,\sigma,\mathcal{L})$. Observe that we have~$s\mid n$. The following proposition then describes the exact value of~$s$, which depends on the behaviour of the normal element~$g$ which lives in degree~3 \cite[theorem~6.8]{MR1086882}.

\begin{proposition}
  \label{proposition:fat-multiplicity}
  If~$(E,\sigma,\mathcal{L})$ is the regular triple associated to an Artin--Schelter regular algebra~$A$ which is finite over its center, then the multiplicity of the fat point modules of~$A$ is
  \begin{equation}
    \label{equation:fat-multiplicity}
    s=
    \begin{cases}
      n & \gcd(n,3)=1 \\
      n/3 & \gcd(n,3)=3.
    \end{cases}
  \end{equation}

  \begin{proof}
    In~\cite[theorem~7.3]{MR1128218} it is shown that~$A[g^{-1}]$ is Azumaya of degree~$s$. Now by~\cite[lemma~5.5.5(i)]{artin-dejong} we have that~$s$ is the order of the automorphism~$\eta$ introduced in~\cite[\S5]{MR1128218}. By~\cite[theorem~5.3.6]{artin-dejong} we have that~$\eta=\sigma^3$, hence~$s$ is~$n$ or~$n/3$ depending on~$\gcd(n,3)$.

    Moreover, by~\cite[lemma~5.5.5(ii)]{artin-dejong} we have that all fat point modules are of multiplicity~$s$.
  \end{proof}
\end{proposition}

\begin{remark}
  It can be shown that a Sklyanin algebra associated to a translation of order~3 has the property that~$\qgr A\cong\coh\mathbb{P}^2$. This case is not considered in the remainder of this paper.
\end{remark}

We will also need the following two facts about fat point modules.
\begin{proposition}
  \label{proposition:fat-g-torsionfree}
  Fat point modules are~$g$\dash torsion free.

  \begin{proof}
    Let~$M$ be a simple graded~$A$\dash module. Then there exists an~$n \in \mathbb{Z}$ such that~$M_i =0$ for all~$i \neq n$. To see this, note that if~$M_i \neq 0$ and~$M_n \neq 0$ for some~$i > n$, then the truncation~$M_{\geq n+1}$ is a non-trivial submodule of~$M$.

    In particular, let~$F$ be a fat point module and~$M \subset F$ a simple graded submodule. By the above and \cref{definition:fat-point-module}(\ref{enumerate:fat-point-module-2}). $M$ is finite dimensional, implying~$M=0$ by \cref{definition:fat-point-module}(\ref{enumerate:fat-point-module-3}).

    Hence~$F$ has a trivial socle. As such we can apply \cite[proposition~7.7(ii)]{MR1128218} from which the lemma follows because~$F$ cannot be an extension of point modules as we assumed~$\pi F \in \qgr A$ to be simple.
  \end{proof}
\end{proposition}

\begin{lemma}
  \label{lemma:triple-shift}
The fat point module~$F$ is invariant under triple degree shifting, i.e.~there exists an isomorphism in~$\gr A$:
  \begin{equation}
    \label{equation:triple-shift}
    F\cong\left(F(3)\right)_{\geq 0}.
  \end{equation}

  \begin{proof}
    This is direct corollary of \cref{proposition:fat-g-torsionfree}: the isomorphism is given by multiplication by~$g$.
  \end{proof}
\end{lemma}

\section{Construction}
\label{section:construction}

\subsection{Noncommutative planes finite over their center}
Consider a 3\dash dimensional quadratic Artin--Schelter-regular algebra~$A$ which is finite over its center~$\ZZ(A)$. Using \cref{theorem:finite-over-center} we know that this is the case precisely when the automorphism in the associated elliptic triple is of finite order.

In this case we can consider
\begin{equation}
  X\coloneqq\Proj\ZZ(A),
\end{equation}
and the sheafification~$\mathcal{R}$ of~$A$ over~$X$. This is a sheaf of noncommutative~$\mathcal{O}_X$\dash algebras, coherent as~$\mathcal{O}_X$\dash module. Often, but not always, we have that~$X\cong\mathbb{P}^2$~\cite[theorem~5.2]{MR1144023}. It is possible to improve this situation by considering a finite cover of~$X$.

The center~$\ZZ(\mathcal{R})$, which is not necessarily~$\mathcal{O}_X$, is a (coherent) sheaf of commutative~$\mathcal{O}_X$\dash algebras, hence we can consider
\begin{equation}
  f\colon Y\coloneqq\relSpec_X\ZZ(\mathcal{R})\to X.
\end{equation}
Because~$\ZZ(\mathcal{R})$ is coherent as an~$\mathcal{O}_X$\dash module the projection map~$f$ is finite.

The main result about~$Y$, for any Artin--Schelter regular algebra finite over its center, is that~$Y$ is isomorphic to~$\mathbb{P}^2$. This was proven:
\begin{enumerate}
  \item by Artin for Sklyanin algebras associated to points of order coprime to~3, where~$X\cong Y$, as mentioned before,
  \item by Smith--Tate for all Sklyanin algebras~\cite{MR1273835},
  \item by Mori for algebras of type~$\mathrm{S}_1$~\cite{MR1624000} (these have a triangle of~$\mathbb{P}^1$'s as their point scheme),
  \item and finally by Van Gastel in complete generality~\cite{MR1880659}, with an analogous proof in~\cite[theorem~5.3.7]{artin-dejong}.
\end{enumerate}

We will denote by~$\Sscr$ the sheaf of algebras on~$Y$ induced by~$\Rscr$, so the situation is described as follows.
\begin{equation}
  \label{equation:central-Proj}
  \begin{tikzcd}
    \mathcal{S} \arrow[d, dashed, dash] & \mathcal{R}=f_*(\mathcal{S}) \arrow[d, dashed, dash] \\
    \mathbb{P}^2\cong Y\coloneqq\relSpec_X\ZZ(\mathcal{R}) \arrow[r, "f"] & X.
  \end{tikzcd}
\end{equation}

The sheaf of algebras~$\mathcal{S}$ has many pleasant properties and will be used in the construction.
\begin{lemma}
  $\Sscr$ is a sheaf of maximal orders on~$\mathbb{P}^2$ of rank~$s^2$, with~$s$ as in \cref{proposition:fat-multiplicity}.

  \begin{proof}
    By the discussion above we have~$Y \cong \mathbb{P}^2$, so that~$\Sscr$ is a maximal order follows from~\cite[proposition~1]{MR1356364}. Observe that the notation in the statement of loc.~cit.~is somewhat unfortunate, and should be taken as in \eqref{equation:central-Proj}.

    It is locally free because it is a reflexive sheaf over a regular scheme of dimension~2, and the statement on the rank follows from \cite[theorem~7.3]{MR1128218}.
  \end{proof}
\end{lemma}
Using this we can decompose~$\mathbb{P}^2$ into a \emph{ramification divisor}~$C$ and its complement, the \emph{Azumaya locus}.

\begin{remark}
  \label{remark:artin-mumford-sequence}
  It is also possible to classify the curves that can appear as ramification divisors for a maximal order on~$\mathbb{P}^2$ using the Artin--Mumford sequence \cite{MR0321934}, as explained in \cite[lemma~1.1(2)]{MR1601190}. This gives the same result as \cref{remark:finite-order-automorphisms}, taking care of the distinction between the point scheme and the ramification divisor.
\end{remark}

The noncommutative plane associated to the algebra~$A$ is the abelian category~$\qgr A$. In the case where~$A$ is finite over its center we have a second interpretation for this category, namely as the category of coherent~$\mathcal{R}$\dash and~$\mathcal{S}$\dash modules.
\begin{lemma}
  \label{lemma:equivalence-qgr-coh}
  There are equivalences of categories
  \begin{equation}
    \qgr A\cong\coh\mathcal{R}\cong\coh\mathcal{S}.
  \end{equation}

  \begin{proof}
    The equivalence~$\QGr A\cong\Qcoh\Rscr$ is given by the restriction of the equivalence~$\widetilde{(-)}: \QGr \ZZ(A)\rightarrow\Qcoh X$. Similarly there is an equivalence of categories~$\Qcoh \ZZ(\Rscr)\cong\Qcoh Y$, and one easily checks that this restricts to~$\Qcoh \Rscr\cong\Qcoh \Sscr$ (see for example~\cite[proposition~3.5]{MR3695056}). This equivalence also restricts to noetherian objects.
  \end{proof}
\end{lemma}

\subsection{Description of the exceptional sequence}
We will use the notation~$\Sscr_i\in\coh\Sscr$ for the images of~$\pi A(i)\in\qgr A$ under the above equivalences. Similarly we fix a fat point module~$F$ and let~$\Fscr \in\coh\Sscr$ be its image. The collection
\begin{equation}
  (\Sscr_0,\Sscr_1,\Sscr_2,\Fscr)
\end{equation}
in~$\derived^\bounded(\Sscr)$ is the noncommutative analogue of~$\Oscr_{\mathbb{P}^2},\Oscr_{\mathbb{P}^2}(1),\Oscr_{\mathbb{P}^2},k(x)$, where~$k(x)$ is the skyscraper in a closed point~$x$. This is not an exceptional collection for~$\derived^\bounded(\mathbb{P}^2)$: we have that~$\Ext^2(k(x),\mathcal{O}_{\mathbb{P}^2}(i))\neq 0$ for all~$i$, but it will become one after blowing up at~$p$.

The point we wish to blow up is the support of the fat point module, considered as an object in~$\coh\Sscr$. This corresponds precisely with a point of~$\mathbb{P}^2\setminus C$, where~$C$ is the ramification divisor of~$\Sscr$. This will give us a new exceptional object, with the appropriate number of morphisms towards it.

We can perform the analogous construction in the noncommutative situation. Let~$x\in\mathbb{P}^2\setminus C$ be the unique closed point in the support of~$\mathcal{F}$, where~$C$ is the ramification locus of~$\mathcal{S}$.

Consider the blowup square
\begin{equation}
  \begin{tikzcd}
    E=\mathbb{P}^1 \arrow[d, "q"] \arrow[r, "j"] & Z=\mathbb{F}_1 \arrow[d, "p"] \\
    x \arrow[r, "i"] & Y=\mathbb{P}^2.
  \end{tikzcd}
\end{equation}
As in \cref{subsection:orlov-for-orders} we will use the notation~$p_\Sscr^*$ for the inverse image functor obtained from the morphism of ringed spaces~$(\mathbb{F}_1,p^*(\Sscr))\to(\mathbb{P}^2,\Sscr)$.

As explained in \cref{example:noncommutative-beilinson} the structure of~$\derived^\bounded(\Sscr)$ is obtained by changing the relations in the quiver according to the generators and relations for the Artin--Schelter regular algebra. This is a well known result, and a noncommutative analogue of Serre's description of the sheaf cohomology of~$\Oscr_{\mathbb{P}^n}(i)$.

\begin{lemma}
  \label{lemma:structure-DbS}
  Let~$A$ be a quadratic~3\dash dimensional Artin--Schelter regular algebra. Then there is a full and strong exceptional collection
  \begin{equation}
    \derived^\bounded(\qgr A)=\left\langle \pi A,\pi A(1),\pi A(2) \right\rangle,
  \end{equation}
  such that
  \begin{equation}
    \begin{aligned}
      \Hom_{\derived^\bounded(\qgr A)}(\pi A(i),\pi A(i+1))&\cong A_1 \\
      \Hom_{\derived^\bounded(\qgr A)}(\pi A(i),\pi A(i+2))&\cong A_2
    \end{aligned}
  \end{equation}
  and the composition law in the quiver is given by the multiplication law~$A_1\otimes_kA_1\twoheadrightarrow A_2$.

  \begin{proof}
    By~\cite[theorem~8.1]{MR1304753} we have
    \begin{equation}
      \begin{aligned}
        \Ext_{\qgr A}^m(\pi A,\pi A(j-i))\cong
        \begin{cases}
          A_{j-i} & m=0 \\
          A_{i-j-3}^\vee & m=2 \\
          0 & m\neq 0,2
        \end{cases}
      \end{aligned}.
    \end{equation}
    Moreover by assumption the algebra is generated in degree~1, hence we have a complete description of the structure of the exceptional collection. That it is full is proven in~\cite[theorem~7.1]{1411.7770v1}.
  \end{proof}
\end{lemma}

Using \cref{lemma:equivalence-qgr-coh} this gives a description for the derived category~$\Sscr$ as
\begin{equation}
  \derived^\bounded(\Sscr)=\left\langle \Sscr_0,\Sscr_1,\Sscr_2 \right\rangle.
\end{equation}

\begin{lemma}
  \label{lemma:Hom-Si-F}
  Let~$F$ be a normalised fat point module for the algebra~$A$ and let~$s$ be as in \cref{proposition:fat-multiplicity}. Then
  \begin{equation}
    \dim_k \left(\Hom_{\qgr A}(\pi A(j),\pi F) \right) = s,
  \end{equation}
  and
  \begin{equation}
    \Ext_{\qgr A}^k(\pi A(j),\pi F)=0
  \end{equation}
  for~$k\geq 1$.

  \begin{proof}
    Using the identities
    \begin{equation}
      \begin{aligned}
        \Hom_{\qgr A}(\pi A(j),\pi F)
        &\cong\Hom_{\qgr A}(\pi A,\pi F(-j)) \\
        F(-j)
        &\cong F(-j+3) & \textrm{ (see \cref{lemma:triple-shift})}
      \end{aligned}
    \end{equation}
    we can assume, without loss of generality, that~$j \leq 0$.

    Recall that
    \begin{equation}
      \Hom_{\qgr A}(\pi A,\pi F(-j))\cong \lim_{i\rightarrow\infty}\Hom_{\gr A }(A_{\geq i}, F(-j)).
    \end{equation}
    We now claim
     \begin{equation}
      \lim_{i\rightarrow\infty}\Hom_{\gr A }(A_{\geq i}, F(-j)) \cong \Hom_{\gr A }(A, F(-j)) \cong F_{-j}.
    \end{equation}
    The lemma follows by combining this claim with \cref{proposition:fat-multiplicity}.

    To prove the claim, note that one can compute the limit by restricting to the directed subsystem~$3\mathbb{N}$. It then suffices to prove that the natural map
    \begin{equation}
      \alpha_i\colon\Hom_{\gr A}(A_{\geq 3i}, F(-j)) \rightarrow \Hom_{\gr A}(A_{\geq 3i+3}, F(-j))
    \end{equation}
    is an isomorphism for all~$i \geq 0$. But this follows as its inverse is given by
    \begin{equation}
      \beta_i\colon\Hom_{\gr A}(A_{\geq 3i+3}, F(-j)) \rightarrow \Hom_{\gr A}(A_{\geq 3i}, F(-j)):\beta_i(\varphi)(x) = g^{-1} \varphi(gx)
    \end{equation}
    where~$g^{-1}\colon F_n\to F_{n-3}$ is the inverse of the isomorphism as in \eqref{equation:triple-shift}.

    To see that the Ext vanish, we use that~$F$ (resp.~$F(-j)$) is finitely generated (by the degree zero part) and has Gelfand--Kirillov dimension 1. By \cite[theorem~4.1]{MR1128218} we can conclude
    \begin{equation}
      \Ext^i_{\gr A}(F,A(j)) = 0 \text{ for all~$j$ and } i=0,1
    \end{equation}
    and~\cite[theorem 8.1]{MR1304753} implies
    \begin{equation}
      \Ext^i_{\qgr A}(\pi F,\pi A(j)) = \Ext^i_{\gr A}(F,A(j)) = 0 \text{ for all~$j$ and } i=0,1
    \end{equation}
    It was proven in~\cite[theorem~2.9.1]{MR2058456} that~$\qgr A$ satisfies the following version of noncommutative Serre duality:
    \begin{equation}
      \Ext^i_{\qgr A}(\pi F, \pi A(j)) \cong \Ext_{\qgr A}^{2-i}(\pi A(j),\pi F(-3))^\vee.
    \end{equation}
    But the latter is isomorphic to~$\Ext_{\qgr A}^{2-i}(\pi A(j),\pi F)^\vee$ by \cref{lemma:triple-shift}, which implies the vanishing of Ext.
  \end{proof}
\end{lemma}

We will also use the following lemma in checking that the exceptional collection is indeed strong, and of the prescribed form.
\begin{lemma}
  \label{lemma:Rp-p-F}
  There exists an isomorphism
  \begin{equation}
    \label{eq:Rp-p-F}
    \RRR p_*\circ p^*(\mathcal{F})\cong\mathcal{F}.
  \end{equation}

  \begin{proof}
    Consider the divisor short exact sequence
    \begin{equation}
      0\to\mathcal{O}_{\mathbb{F}_1}(-E)\to\mathcal{O}_{\mathbb{F}_1}\to j_*(\mathcal{O}_E)\to 0.
    \end{equation}
    Applying the exact functor~$p^*(\mathcal{S})\otimes_{\mathcal{O}_{\mathbb{F}_1}}-$ to it we get a short exact sequene of left~$p^*(\mathcal{S})$\dash modules
    \begin{equation}
      \label{eq:divisor-short-exact-sequence-tensored}
      0
      \to p^*(\mathcal{S})\otimes_{\mathcal{O}_{\mathbb{F}_1}}\mathcal{O}_{\mathbb{F}_1}(-E)
      \to p^*(\mathcal{S})
      \to p^*(\mathcal{S})\otimes_{\mathcal{O}_{\mathbb{F}_1}}j_*(\mathcal{O}_E)
      \to 0.
    \end{equation}

    We have the chain of isomorphisms
    \begin{equation}
      \label{eq:S-restricted-is-nF}
      \begin{aligned}
        &p^*(\mathcal{S})\otimes_{\mathcal{O}_{\mathbb{F}_1}}j_*(\mathcal{O}_E) \\
        &\quad\cong j_*\circ j^*\circ p^*(\mathcal{S}) & \text{projection formula} \\
        &\quad\cong j_*\circ q^*\circ i^*(\mathcal{S}) & \text{functoriality} \\
        &\quad\cong p^*\circ i_*\circ i^*(\mathcal{S}) & \text{base change for affine morphisms} \\
        &\quad\cong p^*(\mathcal{F}^{\oplus s})
      \end{aligned}
    \end{equation}
    where the last step uses that~$i^*(\mathcal{S})\cong\Mat_s(k)$ is the direct sum of the~$s$\dash dimensional representation corresponding to the fat point~$F$.

    Because the first two terms in \eqref{eq:divisor-short-exact-sequence-tensored} are~$p_*$\dash acyclic, so is the third and its direct summands~$p^*(\mathcal{F})$. Hence we can use the projection formula
    \begin{equation}
      p_*(p^*(\mathcal{S})\otimes_{p^*(\mathcal{S})}p^*(\mathcal{F}))\cong p_*\circ p^*(\mathcal{S})\otimes_{\mathcal{S}}\mathcal{F}\cong\mathcal{F}.
    \end{equation}
  \end{proof}
\end{lemma}

We can now prove the main theorem of this paper, which gives a construction of noncommutative surfaces with prescribed Grothendieck group. It uses a semiorthogonal decomposition that generalises Orlov's blowup formula, and which is proved in some generality in \cref{subsection:orlov-for-orders}.

\begin{theorem}
  \label{theorem:main-theorem}
  Let~$A$ be a quadratic~3\dash dimensional Artin--Schelter regular algebra, finite over its center. Let~$F$ be a fat point module of~$A$. Let~$p\colon\Bl_x\mathbb{P}^2\to\mathbb{P}^2$ be the blowup in the point~$x$ which is the support of the~$\Sscr$\dash module~$\Fscr$ as a sheaf on~$\mathbb{P}^2$. Let~$s$ be the integer as in \eqref{equation:fat-multiplicity}. Then
  \begin{equation}
    \derived^\bounded(p^*\mathcal{S})=\left\langle \LLL p_\Sscr^*\Sscr_0,\LLL p_\Sscr^*\Sscr_1,\LLL p_\Sscr^*\Sscr_2,p_\Sscr^*\Fscr \right\rangle
  \end{equation}
  is a full and strong exceptional collection, whose Gram matrix is of type \eqref{equation:family-P2-matrix}, where~$m=s$.

  \begin{proof}
    By \cref{theorem:orlov-for-orders} we obtain that the collection is indeed a full and strong exceptional collection, where we use that~$\mathcal{F}$ can be considered as the (noncommutative) skyscraper sheaf for~$\Sscr$, because we are in the Azumaya locus.

    The structure of~$\langle\LLL p_\Sscr^*\Sscr_0,\LLL p_\Sscr^*\Sscr_1,\LLL p_\Sscr^*\Sscr_2\rangle$ is described in \cref{lemma:structure-DbS} using the fully faithfulness of~$\LLL p_\Sscr^*$ from \cref{lemma:Lp-fully-faithful}.

    Finally using \cref{lemma:Rp-p-F} we get that
    \begin{equation}
      \begin{aligned}
        \Hom_{\derived^\bounded(p^*\Sscr)}(\LLL p_\Sscr^*\Sscr_i,p^*\Fscr)
        &\cong\Hom_{\derived^\bounded(\Sscr)}(\Sscr_i,\Fscr) \\
        &\cong\Hom_{\Sscr}(\Sscr_i,\Fscr)
      \end{aligned}
    \end{equation}
    which is~$s$\dash dimensional by \cref{lemma:Hom-Si-F}, and similarly we get that there are no forward Ext's to conclude that the collection is indeed strong.
  \end{proof}
\end{theorem}

\begin{remark}
  There are three degrees of freedom in this construction: generically the point scheme is an elliptic curve for which we have the~$j$\dash line as moduli space, for each curve there are only finitely many torsion automorphisms, and then there is the choice of a point in~$\mathbb{P}^2\setminus C$. There are only finitely many automorphisms of~$\mathbb{P}^2$ that preserve~$C$, so we get three degrees of freedom. This is the expected number, using \cite{1705.06098v1}, where a formula for~$\dim_k\mathrm{HH}^2-\dim_k\mathrm{HH}^1$ is given in terms of the number of exceptional objects.
\end{remark}

\begin{remark}
  The derived category~$\derived^\bounded(p^*\Sscr)$ comes with its standard t-structure. By \cite{MR738217} we have that Serre duality takes on the form that it does for ordinary smooth and projective varieties. This means that the Serre functor is compatible with the t-structure, as required in Bondal's definition of geometric t-structure~\cite{mpim-93-67}. In particular, it is an example of a noncommutative variety of dimension~2\ in this sense.
\end{remark}

\begin{remark}
  Using \cite{MR3545926} and the full and strong exceptional collection from \cref{theorem:main-theorem} there exists an embedding of~$\derived^\bounded(p^*\mathcal{S})$ into the derived category of a smooth projective variety. Now in the spirit of \cite{MR3488782,1605.02795v1,MR3397451,1612.02241v1,MR3713871} it is an interesting question whether there exists a natural embedding, i.e.~where the smooth projective variety is associated to~$p^*\Sscr$ in a natural way. An obvious candidate would be the Brauer--Severi scheme of the maximal order, and indeed in the case where the automorphism is of order~2 there exists a fully faithful embedding into~$\derived^\bounded(\BS(p^*\Sscr))$ by \cite[\S6]{MR2905553} and \cite{MR2419925}, as the maximal order is the even part of a sheaf of Clifford algebras. What happens for the more general case and the study of the derived category of the Brauer--Severi scheme in this situation, is left for future work.
\end{remark}

\subsection{Orlov's blowup formula for orders}
\label{subsection:orlov-for-orders}
The main ingredient in the construction of \cref{theorem:main-theorem} is the observation that it is possible to generalise Orlov's blowup formula \cite[theorem~4.3]{MR1208153} to a sufficiently nice noncommutative setting where we blow up a point on the underlying variety, and pull back the sheaf of algebras to the blown up variety. This is a result of independent interest.

\begin{theorem}
  \label{theorem:orlov-for-orders}
  Let~$X$ be a smooth quasiprojective variety. Let~$\mathcal{A}$ be a locally free sheaf of orders of degree~$n$ on~$X$ such that~$\gldim\mathcal{A}<+\infty$. Let~$Y$ be a smooth closed subvariety such that it does not meet the ramification locus of~$\mathcal{A}$. Assume moreover that~$\mathcal{A}|_Y\cong\Mat_n(\mathcal{O}_Y)$. Consider the blowup square
  \begin{equation}
    \begin{tikzcd}
      E\coloneqq\mathbb{P}(\mathrm{N}_YX) \arrow[d, "q"] \arrow[r, "j"] & Z\coloneqq \Bl_YX \arrow[d, "p"] \\
      Y \arrow[r, "i"] & X
    \end{tikzcd}
  \end{equation}
  and its noncommutative analogue obtained by pulling back the sheaf of algebras~$\mathcal{A}$
  \begin{equation}
    \begin{tikzcd}
      (E,\Mat_n(\mathcal{O}_E)) \arrow[d, "q_{\mathcal{A}}"] \arrow[r, "j_{\mathcal{A}}"] & (Z,p^*\mathcal{A}) \arrow[d, "p_{\mathcal{A}}"] \\
      (Y,\Mat_n(\mathcal{O}_Y)) \arrow[r, "i_{\mathcal{A}}"] & (X,\mathcal{A})
    \end{tikzcd}.
  \end{equation}

  Then we have a semiorthogonal decomposition
  \begin{equation}
    \label{equation:nc-blowup-sod}
    \derived^\bounded(Z,p^*(\mathcal{A}))
    =
    \left\langle
      \derived^\bounded(X,\mathcal{A}),
      \derived^\bounded(Y),
      \ldots,
      \derived^\bounded(Y)
    \right\rangle
  \end{equation}
  where the first component is embedded using the functor~$\LLL p_{\Ascr}^*$, and the subsequent components by~$j_{\Ascr,*}(q_{\Ascr}^*(-)\otimes\mathcal{O}_E(kE))$, for~$k=0,\ldots,\codim_XY-2$.
\end{theorem}

\begin{remark}
  In this case~$p^*\mathcal{A}$ is automatically of finite global dimension.
\end{remark}

\begin{remark}
  The case where~$\mathcal{A}=\mathcal{O}_X$ is Orlov's blowup formula. Observe that his proof works verbatim for a smooth quasiprojective variety as all morphisms are projective, so the bounded derived category is preserved throughout. We will use this in the proof of \cref{theorem:orlov-for-orders}.
\end{remark}

We can prove this result by bootstrapping the original proof. To do so we will need generalisations of some standard results in algebraic geometry such as the adjunction between (derived) pullback and direct image, or the projection formula. A reference for these in the setting of Azumaya algebras can be found in~\cite[\S10]{MR2238172}. We will only need results that do not depend on the algebras being Azumaya, hence in the words of remark~10.5 of op.~cit.~we are working with noncommutative finite flat (and not \'etale) coverings.

\begin{lemma}
  \label{lemma:Lp-fully-faithful}
  The functor~$\LLL p_{\mathcal{A}}^*$ is fully faithful.

  \begin{proof}
    In the commutative setting this is proven using the derived projection formula and the fact that~$\RRR p_*\circ\LLL p^*(\mathcal{O}_X)\cong\mathcal{O}_X$. In the noncommutative setting the appropriate projection formula is given as the first isomorphism in~\cite[lemma~10.12]{MR2238172} taking into account that we have a bimodule structure, whilst the isomorphism
    \begin{equation}
      \RRR p_{\mathcal{A},*}\circ\LLL p_{\mathcal{A}}^*(\mathcal{A})\cong\mathcal{A}
    \end{equation}
    follows from the third isomorphism in loc.~cit.
  \end{proof}
\end{lemma}

By the assumption that~$\mathcal{A}|_Y\cong\Mat_n(\mathcal{O}_Y)$ and the projection formula as given in the third isomorphism of loc.~cit.~we have that~\cite[lemma~4.2]{MR1208153} goes through as stated. In particular we only need to check that the semiorthogonal decomposition is indeed full.

\begin{proof}[of \cref{theorem:orlov-for-orders}]
  We can mimick the proof of the first part of~\cite[theorem~4.3]{MR1208153}. Consider an object in the right orthogonal of \eqref{equation:nc-blowup-sod}, in particular it is right orthogonal to~$\LLL p_{\mathcal{A}}^*(\derived^\bounded(X,\mathcal{A}))$. As Serre duality for sheaves of maximal orders takes on the expected form by~\cite[corollary~2]{MR738217} we get that~$\RRR p_*$ of this object is indeed zero, and therefore that it is contained in the minimal full subcategory containing the image of~$\derived(E,\Mat_n(\mathcal{O}_E))$.

  Now use that blowups commute with flat base change, in particular we can take an \'etale neighbourhood of the exceptional divisor that splits~$\mathcal{A}$ and such that its image in~$X$ and the ramification divisor are disjoint. Then we are in the usual setting of Orlov's blowup formula (up to Morita equivalence), and we can use the usual proof to conclude that the object is indeed zero.
\end{proof}

Two remarks are in order, which are already important in the case of blowing up a point on a surface.
\begin{remark}
  If we were to blow up a point \emph{on the ramification divisor}, then the resulting algebra is not necessarily of finite global dimension. Considering the case of a Sklyanin algebra associated to a point of order~2, we have that the complete local structure of this algebra at the point on the intersection of the exceptional divisor and the ramification locus is given by
  \begin{equation}
    \begin{pmatrix}
      R & R \\
      (xy) & R
    \end{pmatrix}
  \end{equation}
  where~$R=k[[x,y]]$. One then checks that the module~$\begin{psmallmatrix} 0 \\ R/(x) \end{psmallmatrix}$ has a periodic minimal projective resolution of the form
  \begin{equation}
    \ldots
    \rightarrow \begin{pmatrix} (xy) \\ (xy) \end{pmatrix} \oplus \begin{pmatrix} (x) \\ (x^2y) \end{pmatrix}
    \overset{\psi'}{\rightarrow} \begin{pmatrix} (x) \\ (x) \end{pmatrix} \oplus \begin{pmatrix} R \\ (xy) \end{pmatrix}
    \overset{\varphi}{\rightarrow} \begin{pmatrix} R \\ R \end{pmatrix}
    \overset{\psi}{\rightarrow} \begin{pmatrix} 0 \\ R/(x) \end{pmatrix}
    \rightarrow 0
  \end{equation}
  To see that this resolution is in fact periodic, note that~$\Ker(\psi') \cong \Ker(\psi)$ as
  \begin{equation}
    \Ker(\psi) = \begin{pmatrix} R \\ (x) \end{pmatrix} \textrm{ and } \Ker(\psi') \cong \begin{pmatrix} (xy) \\ (x^2y) \end{pmatrix} =  xy \begin{pmatrix} R \\ (x) \end{pmatrix}.
  \end{equation}

  This description also shows that the result is no longer a maximal order, and there is a choice of embedding, as explained in~\cite[\S4]{MR2180454}. Without the embedding in a maximal order one does not expect a meaningful semiorthogonal decomposition. In this paper we do not need to take a maximal order containing the pullback as it is already maximal using the Auslander--Goldman criterion.
\end{remark}

\begin{remark}
  In the construction of \cite{MR1846352} a point \emph{on the point scheme} is blown up. For an algebra finite over its center this is not the same as the ramification curve, these two curves are only isogeneous. In the context of the previous remark, the difference is measured by the choice of a maximal order containing the pullback.
\end{remark}

\section{Properties of the maximal orders}
\label{section:order-properties}
In the minimal model program for orders on surfaces as studied in~\cite{MR2180454,artin-dejong} there exists the notion of del Pezzo orders, and (half-)ruled orders. The orders we have constructed are obviously not minimal, but they give rise to interesting examples in the study of maximal orders.

In the commutative case the surface~$\Bl_x\mathbb{P}^2=\mathbb{F}_1$ is both del Pezzo and ruled. We are considering orders \emph{on} this surface, and for the value of~$m=2$ in the classification we obtain that it is both \emph{del Pezzo} and \emph{half ruled}, as explained in \cref{proposition:m-2-del-pezzo} and \cref{proposition:m-2-half-ruled}.

This latter notion is introduced by Artin, to describe a class of orders which is not ruled, but whose cohomological properties mimick those of ruled orders.

\subsection{The case \texorpdfstring{$m=2$}{m=2} is del Pezzo}
In this section we quickly recall the notion of del Pezzo order, and show that the intuition from numerical blowups from \cite{1607.04246v1} agrees with the a priori independent notion of del Pezzo order, introduced in~\cite[\S3]{MR1954458}.

Throughout we let~$\mathcal{A}$ be a maximal order on a smooth projective surface~$S$.
\begin{definition}
  The \emph{canonical sheaf of~$\mathcal{A}$} is the~$\mathcal{A}$-bimodule
  \begin{equation}
    \omega_{\mathcal{A}}:=\mathcal{H}\mathrm{om}_{\mathcal{O}_S}(\mathcal{A},\omega_S).
  \end{equation}
\end{definition}
Now denote~$\omega_{\mathcal{A}}^*:=\mathcal{H}\mathrm{om}_{\mathcal{A}}(\omega_{\mathcal{A}},\mathcal{A})$, whereas~$\mathcal{F}^\vee$ is used for~$\mathcal{H}\mathrm{om}_{\mathcal{O}_S}(\mathcal{F},\mathcal{O}_S)$, hence the reflexive hull is denoted~$\mathcal{F}^{\vee\vee}$.

\begin{definition}
  Let~$\mathcal{L}$ be an invertible~$\mathcal{A}$\dash bimodule, which is moreover~$\mathbb{Q}$\dash Cartier, i.e.~$(\mathcal{L}^{\otimes n})^{\vee\vee}$ is again invertible for some~$n$. Then~$\mathcal{L}$ is \emph{ample} if
  \begin{equation}
    \mathrm{R}^q\Hom_{\mathcal{A}}(\mathcal{A},(\mathcal{L}^{\otimes k})^{\vee\vee}\otimes\mathcal{F})\cong\mathrm{H}^q(S,(\mathcal{L}^{\otimes k})^{\vee\vee}\otimes\mathcal{F})
  \end{equation}
  is zero for~$q\geq 1$ and~$k\gg0$, where the isomorphism is induced by applying the forgetful functor.
\end{definition}

Then analogous to the commutative situation we define
\begin{definition}
  The maximal order~$\mathcal{A}$ is \emph{del Pezzo} if~$\omega_{\mathcal{A}}^\vee$ is ample.
\end{definition}

In particular, by~\cite[lemma~8]{MR1954458} it suffices to check that the divisor
\begin{equation}
  \mathrm{K}_{\mathcal{A}}=\mathrm{K}_S+\sum_{i=1}^n\left( 1-\frac{1}{e_i} \right)C_i
\end{equation}
is anti-ample: the del Pezzoness only depends on the center and the ramification data.

\begin{proposition}
  \label{proposition:m-2-del-pezzo}
  Let~$\mathcal{A}$ be the pullback of a maximal order of degree~$m$ on~$\mathbb{P}^2$ ramified on a cubic curve along the blowup~$\mathbb{F}_1\to\mathbb{P}^2$ in a point outside the ramification locus. Then~$\mathcal{A}$ is del Pezzo if and only if~$m=2$.

  \begin{proof}
    If we denote~$\operatorname{Pic}\mathbb{F}_1=\mathbb{Z}H\oplus\mathbb{Z}E$, such that~$H^2=1$,~$H\cdot E=0$ and~$E^2=-1$, then
    \begin{equation}
      \begin{aligned}
        \mathrm{K}_{\mathcal{A}}
        &=-3 p^*(H)+E+\left( 1-\frac{1}{m} \right)3p^*(H) \\
        &=\frac{-3}{m}p^*(H)+E
      \end{aligned}
    \end{equation}
    because the ramification data for the pullback is the pullback of the ramification data, which is a cyclic cover of degree~$n\geq 2$ of an elliptic curve.

    By the Kleiman criterion for ampleness we need to check that~$-\mathrm{K}_{\mathcal{A}}\cdot C\geq 0$, for~$C$ in the Mori--Kleiman cone of~$\mathbb{F}_1$. This cone is spanned by a fibre~$f$ and the section~$C_0$ of the projection~$\mathbb{F}_1\to\mathbb{P}^1$. We have that~$p^*(H)=C_0+f$ and~$E=C_0$ in the translation between the canonical bases for~$\Pic(\mathbb{F}_1)$ and~$\Pic(\Bl_x\mathbb{P}^2)$.

    Using the description of the intersection form on a ruled surface we obtain
    \begin{equation}
      \begin{aligned}
        -\mathrm{K}_{\mathcal{A}}\cdot f&=\frac{3}{m}-1, \\
        -\mathrm{K}_{\mathcal{A}}\cdot C_0&=1.
      \end{aligned}
    \end{equation}
    The first intersection number is positive if and only if~$m=2$. The second intersection number is always positive.
  \end{proof}
\end{proposition}

\begin{remark}
  The computation for the del Pezzoness of the numerical blowup reduces to the same equation (up to multiplication by~$m^2$).
\end{remark}

\begin{remark}
  For the case~$m=2$ the equations for a Sklyanin algebra from \cref{example:sklyanin} take on a particularly easy form
  \begin{equation}
    \left\{
      \begin{aligned}
        xy+yx+cz^2&=0 \\
        yz+zy+cx^2&=0 \\
        zx+xz+cy^2&=0
      \end{aligned}
    \right.
  \end{equation}
  where~$c^3\neq 0,1,-8$,~\cite[theorem~3.1]{MR3366574}. The case~$c=0$ corresponds to an order associated to a so called skew polynomial algebra, and the ramification curve is a triangle of~$\mathbb{P}^1$'s. If~$c^3=1,-8$ the ramification curve is the union of a conic and a line in general position. We will come back to this situation in \cite{belmans-presotto-vandenbergh-comparison}.
\end{remark}

\subsection{The case \texorpdfstring{$m=2$}{m=2} is half ruled}
\label{subsection:half-ruled}
In the following definition, the curve of genus~$0$ will be a curve over the function field of the base curve of a ruled surface, i.e.~if we consider~$\pi\colon S\to C$ over the field~$k$, then~$K$ will be the function field of~$\mathbb{P}_{k(C)}^1$.
\begin{definition}
  Let~$\mathcal{A}$ be a maximal order in a central simple algebra of degree~2 over the function field~$K$ of a curve~$X$ of genus~$0$. If~$\mathcal{A}$ is ramified in~3 points, with ramification degree~2\ in each point, then we say that~$\mathcal{A}$ is \emph{half ruled}.
\end{definition}

\begin{remark}
  The case where the ramification is of type~$(e,e)$ is the \emph{ruled} case.
\end{remark}

It is shown in~\cite[proposition~4.2.4]{artin-dejong} that being (half-)ruled is equivalent to the Euler characteristic~$\chi(X,\mathcal{A})$ of the coherent sheaf~$\mathcal{A}$ being positive.

The following definition seems to be missing as such from the literature, but it is used implicitly in~\cite{artin-dejong}.
\begin{definition}
  Let~$\mathcal{A}$ be a maximal order on a ruled surface~$S\to C$. Then we say that it is \emph{half ruled} if the fiber of the order over the generic point of~$C$ is half ruled.
\end{definition}

\begin{proposition}
  \label{proposition:m-2-half-ruled}
  The sheaf of maximal orders constructed for the case~$m=2$ is half ruled.

  \begin{proof}
    The ramification divisor on~$\mathbb{P}^2$ being a cubic curve we get that the generic intersection of the fibre of the ruling with the inverse image of the ramification divisor in~$\mathbb{F}_1$ is~3, which proves the claim.
  \end{proof}
\end{proposition}

\subsection{The case \texorpdfstring{$m=3$}{m=3} is elliptic}
We will reuse the notation of \cref{subsection:half-ruled}.
\begin{definition}
  Let~$\mathcal{A}$ be a maximal order in a central simple algebra of degree~3 over the function field~$K$ of a curve~$X$ of genus~$0$. If~$\mathcal{A}$ is ramified in~3 points, with ramification degree~3\ in each point, then we say that~$\mathcal{A}$ is \emph{elliptic}.
\end{definition}

The following proposition is proven in the same way as \cref{proposition:m-2-half-ruled}.
\begin{proposition}
  \label{proposition:m-3-elliptic}
  The sheaf of maximal orders constructed for the case~$m=3$ is elliptic.
\end{proposition}

\newpage

\bibliographystyle{plain}
\bibliography{construction}

\begin{thebibliography}{10}

\bibitem{1411.7770v1}
Tarig Abdelgadir, Shinnosuke Okawa, and Kazushi Ueda.
\newblock Compact moduli of noncommutative projective planes.

\bibitem{artin-dejong}
M.~Artin, D.~Chan, A.~J. de~Jong, and M.~Lieblich.
\newblock Terminal orders on surfaces.
\newblock 2014.

\bibitem{MR1144023}
Michael Artin.
\newblock Geometry of quantum planes.
\newblock In {\em Azumaya algebras, actions, and modules ({B}loomington, {IN},
  1990)}, volume 124 of {\em Contemp. Math.}, pages 1--15. Amer. Math. Soc.,
  Providence, RI, 1992.

\bibitem{MR0321934}
Michael Artin and David Mumford.
\newblock Some elementary examples of unirational varieties which are not
  rational.
\newblock {\em Proc. London Math. Soc. (3)}, 25:75--95, 1972.

\bibitem{MR917738}
Michael Artin and William Schelter.
\newblock Graded algebras of global dimension {$3$}.
\newblock {\em Adv. in Math.}, 66(2):171--216, 1987.

\bibitem{MR1086882}
Michael Artin, John Tate, and Michel Van~den Bergh.
\newblock Some algebras associated to automorphisms of elliptic curves.
\newblock In {\em The {G}rothendieck {F}estschrift, {V}ol.\ {I}}, volume~86 of
  {\em Progr. Math.}, pages 33--85. Birkh\"auser Boston, Boston, MA, 1990.

\bibitem{MR1128218}
Michael Artin, John Tate, and Michel Van~den Bergh.
\newblock Modules over regular algebras of dimension {$3$}.
\newblock {\em Invent. Math.}, 106(2):335--388, 1991.

\bibitem{MR1304753}
Michael Artin and James Zhang.
\newblock Noncommutative projective schemes.
\newblock {\em Adv. Math.}, 109(2):228--287, 1994.

\bibitem{MR509388}
A.~A. Be{\u\i}linson.
\newblock Coherent sheaves on {${\bf P}^{n}$}\ and problems in linear algebra.
\newblock {\em Funktsional. Anal. i Prilozhen.}, 12(3):68--69, 1978.

\bibitem{belmans-presotto-vandenbergh-comparison}
P.~Belmans, D.~Presotto, and M.~Van~den Bergh.
\newblock Comparing two constructions of noncommutative surfaces of rank 4.

\bibitem{1705.06098v1}
Pieter Belmans.
\newblock Hochschild cohomology of noncommutative planes and quadrics.

\bibitem{1605.02795v1}
Pieter Belmans and Theo Raedschelders.
\newblock Noncommutative quadrics and {H}ilbert schemes of points.

\bibitem{mpim-93-67}
A.~Bondal.
\newblock Non-commutative deformations and {P}oisson brackets on projective
  spaces.
\newblock MPIM/93-67, 1993.

\bibitem{MR1230966}
Alexey Bondal and Alexander Polishchuk.
\newblock Homological properties of associative algebras: the method of
  helices.
\newblock {\em Izv. Ross. Akad. Nauk Ser. Mat.}, 57(2):3--50, 1993.

\bibitem{MR3695056}
Igor Burban, Yuriy Drozd, and Volodymyr Gavran.
\newblock Minors and resolutions of non-commutative schemes.
\newblock {\em Eur. J. Math.}, 3(2):311--341, 2017.

\bibitem{MR2180454}
Daniel Chan and Colin Ingalls.
\newblock The minimal model program for orders over surfaces.
\newblock {\em Invent. Math.}, 161(2):427--452, 2005.

\bibitem{MR2905553}
Daniel Chan and Colin Ingalls.
\newblock Conic bundles and {C}lifford algebras.
\newblock In {\em New trends in noncommutative algebra}, volume 562 of {\em
  Contemp. Math.}, pages 53--75. Amer. Math. Soc., Providence, RI, 2012.

\bibitem{MR1954458}
Daniel Chan and Rajesh~S. Kulkarni.
\newblock del {P}ezzo orders on projective surfaces.
\newblock {\em Adv. Math.}, 173(1):144--177, 2003.

\bibitem{MR3366574}
Kevin De~Laet.
\newblock Graded {C}lifford algebras of prime global dimension with an action
  of {$H_p$}.
\newblock {\em Comm. Algebra}, 43(10):4258--4282, 2015.

\bibitem{MR2058456}
Koen de~Naeghel and Michel Van~den Bergh.
\newblock Ideal classes of three-dimensional {S}klyanin algebras.
\newblock {\em J. Algebra}, 276(2):515--551, 2004.

\bibitem{1607.04246v1}
Louis de~Thanhoffer~de V\"olcsey and Michel~Van den Bergh.
\newblock On an analogue of the {M}arkov equation for exceptional collections
  of length 4.

\bibitem{1503.03992v4}
Louis de~Thanhoffer~de V\"olcsey and Dennis Presotto.
\newblock Homological properties of a certain noncommutative del pezzo surface.

\bibitem{1612.02241v1}
Anton Fonarev and Alexander Kuznetsov.
\newblock Derived categories of curves as components of fano manifolds.

\bibitem{MR3397451}
Andreas Krug and Pawel Sosna.
\newblock On the derived category of the {H}ilbert scheme of points on an
  {E}nriques surface.
\newblock {\em Selecta Math. (N.S.)}, 21(4):1339--1360, 2015.

\bibitem{MR1286839}
Sergej Kuleshov and Dmitri Orlov.
\newblock Exceptional sheaves on {D}el {P}ezzo surfaces.
\newblock {\em Izv. Ross. Akad. Nauk Ser. Mat.}, 58(3):53--87, 1994.

\bibitem{MR2238172}
Alexander Kuznetsov.
\newblock Hyperplane sections and derived categories.
\newblock {\em Izv. Ross. Akad. Nauk Ser. Mat.}, 70(3):23--128, 2006.

\bibitem{MR2419925}
Alexander Kuznetsov.
\newblock Derived categories of quadric fibrations and intersections of
  quadrics.
\newblock {\em Adv. Math.}, 218(5):1340--1369, 2008.

\bibitem{MR3691718}
Alexander Kuznetsov.
\newblock Exceptional collections in surface-like categories.
\newblock {\em Mat. Sb.}, 208(9):116--147, 2017.

\bibitem{MR1356364}
Lieven Le~Bruyn.
\newblock Central singularities of quantum spaces.
\newblock {\em J. Algebra}, 177(1):142--153, 1995.

\bibitem{MR1624000}
Izuru Mori.
\newblock The center of some quantum projective planes.
\newblock {\em J. Algebra}, 204(1):15--31, 1998.

\bibitem{MR3713871}
M.~S. Narasimhan.
\newblock Derived categories of moduli spaces of vector bundles on curves.
\newblock {\em J. Geom. Phys.}, 122:53--58, 2017.

\bibitem{MR1208153}
Dmitri Orlov.
\newblock Projective bundles, monoidal transformations, and derived categories
  of coherent sheaves.
\newblock {\em Izv. Ross. Akad. Nauk Ser. Mat.}, 56(4):852--862, 1992.

\bibitem{MR3488782}
Dmitri Orlov.
\newblock Geometric realizations of quiver algebras.
\newblock {\em Proc. Steklov Inst. Math.}, 290(1):70--83, 2015.

\bibitem{MR3545926}
Dmitri Orlov.
\newblock Smooth and proper noncommutative schemes and gluing of {DG}
  categories.
\newblock {\em Adv. Math.}, 302:59--105, 2016.

\bibitem{MR1273836}
S.~Paul Smith.
\newblock The four-dimensional {S}klyanin algebras.
\newblock In {\em Proceedings of {C}onference on {A}lgebraic {G}eometry and
  {R}ing {T}heory in honor of {M}ichael {A}rtin, {P}art {I} ({A}ntwerp, 1992)},
  volume~8, pages 65--80, 1994.

\bibitem{MR1273835}
S.~Paul Smith and John Tate.
\newblock The center of the {$3$}-dimensional and {$4$}-dimensional {S}klyanin
  algebras.
\newblock In {\em Proceedings of {C}onference on {A}lgebraic {G}eometry and
  {R}ing {T}heory in honor of {M}ichael {A}rtin, {P}art {I} ({A}ntwerp, 1992)},
  volume~8, pages 19--63, 1994.

\bibitem{MR1601190}
Michel Van~den Bergh.
\newblock Division algebras on {${\bf P}^2$} of odd index, ramified along a
  smooth elliptic curve are cyclic.
\newblock In {\em Alg\`ebre non commutative, groupes quantiques et invariants
  ({R}eims, 1995)}, volume~2 of {\em S\'emin. Congr.}, pages 43--53. Soc. Math.
  France, Paris, 1997.

\bibitem{MR1846352}
Michel Van~den Bergh.
\newblock Blowing up of non-commutative smooth surfaces.
\newblock {\em Mem. Amer. Math. Soc.}, 154(734):x+140, 2001.

\bibitem{MR2836401}
Michel Van~den Bergh.
\newblock Noncommutative quadrics.
\newblock {\em Int. Math. Res. Not. IMRN}, (17):3983--4026, 2011.

\bibitem{MR738217}
Michel Van~den Bergh and Jan Van~Geel.
\newblock A duality theorem for orders in central simple algebras over function
  fields.
\newblock {\em J. Pure Appl. Algebra}, 31(1-3):227--239, 1984.

\bibitem{MR1880659}
Martine Van~Gastel.
\newblock On the center of the {P}roj of a three dimensional regular algebra.
\newblock {\em Comm. Algebra}, 30(1):1--25, 2002.

\end{thebibliography}

\end{document}